\documentclass{article}
\usepackage{blindtext}
\usepackage[utf8]{inputenc}

\usepackage{graphics}
\usepackage[
backend=biber,
style=ieee,
sorting=nyt
]{biblatex}
\usepackage{graphicx}
\usepackage{mathtools}
\usepackage{amssymb}
\usepackage{amsfonts}
\usepackage{amsthm}
\usepackage{amsmath}
\usepackage{txfonts}
\usepackage{float}
\usepackage{caption}
\usepackage{subcaption}
\usepackage{comment}
\usepackage{enumitem}
\usepackage[mathlines]{lineno}
\usepackage{tikz-cd}

\usepackage{setspace}

\usepackage{hyperref}
\hypersetup{
	colorlinks=true,
	linkcolor=blue,
	filecolor=magenta,      
	urlcolor=cyan,
	pdfpagemode=FullScreen,
}

\usepackage{cases}
\usepackage{orcidlink}

\DeclareMathOperator{\diag}{diag}
\DeclareMathOperator{\Sym}{Sym}

\DeclareMathOperator{\Gr}{Gr}
\DeclareMathOperator{\Ker}{Ker}
\DeclareMathOperator{\rk}{rk}
\DeclareMathOperator{\Mat}{Mat}
\DeclareMathOperator{\Cov}{Cov}

\DeclareMathOperator{\Span}{Span}

\DeclareMathOperator{\Imm}{Im}
\DeclareMathOperator{\ver}{ver}
\DeclareMathOperator{\hor}{hor}
\DeclareMathOperator{\St}{St}
\DeclareMathOperator{\Skew}{Skew}
\DeclareMathOperator{\Hom}{Hom}
\DeclareMathOperator{\pr}{pr}

\newcommand{\cl}[1]{\mathcal{#1}}

\newcommand{\Tr}{\textup{\textrm{Tr }}}

\newtheorem{proposition}{Proposition}
\newtheorem{theorem}{Theorem}
\newtheorem{corollary}{Corollary}
\newtheorem{lemma}{Lemma}
\theoremstyle{definition}
\newtheorem{example}{Example}
\theoremstyle{definition}
\newtheorem{definition}{Definition}
\theoremstyle{remark}
\newtheorem*{remark}{Remark}

\def\R{\mathbb{R}}
\def\P{\mathbb{P}}

\allowdisplaybreaks

\title{An associated bundle approach to the Bures--Wasserstein geometry of fixed rank covariance matrices}
\author{Leonardo Marconi
	\thanks{
		Dipartimento di Scienze Statistiche Paolo
		Fortunati, Università di Bologna, Bologna, Italy. \textbf{e-mail}:
		leonardo.marconi5@unibo.it}
	\orcidlink{0009-0000-1597-7875}
}
\date{\today}

\allowdisplaybreaks

\begin{document}
	
	\maketitle
	\bigskip
	
\begin{abstract}
	The Bures--Wasserstein geometry of covariance matrices provides a canonical distance on the statistical manifold of centred Gaussian measures and lies at the intersection of information geometry, quantum information, and optimal transport. The space of covariance matrices admits a natural stratified structure whose strata consist of fixed-rank covariance matrices. In this paper we focus on the rank-$k$ stratum $\Sym^+(n,k)$ and revisit its geometry through the diffeomorphic associated-bundle model $\Sym^+(n,k)\cong\St(n,k)\times_{O(k)}\Sym^{+}(k)$. Working in this bundle picture, we (i) derive a system of differential equations for Bures--Wasserstein geodesics, (ii) prove that the fibers are totally geodesic and (iii) establish a one-to-one correspondence between Grassmannian logarithms and Bures--Wasserstein logarithms on $\Sym^+(n,k)$, and hence between minimizing geodesics in the two spaces. This alternative viewpoint clarifies the role of the underlying base $\Gr(k,n)$ in the Bures--Wasserstein geometry of low-rank covariance matrices and sets the stage for further investigations into structured covariance models.
\end{abstract}

	\textbf{Keywords}: fixed rank covariance matrices, Bures--Wasserstein, Grassmannian,
	geodesics, fiber bundle.

\section{Introduction}
Covariance matrices arise in many areas of statistics, physics, and engineering, as they encode the second-order structure of multivariate data.  
In particular, the intrinsic geometry of covariance matrices under the Bures--Wasserstein metric is of great relevance in information geometry, quantum information (\cite{bures1969extension}, \cite{helstrom1967minimum}) and optimal transport (\cite{dowson1982frechet}, \cite{olkin1982distance}). A deeper geometric understanding of the space of covariance matrices translates into stronger theoretical guarantees and more efficient algorithms across applications such as covariance estimation, low-rank approximation, and state-space modelling. Recent statistical literature has explicitly exploited this geometric framework to develop estimation procedures and regularization schemes for covariance matrices, see e.g.\ \cite{chandrasekaran2010latent,bernardi2022log}.

The smooth manifold of fixed-rank covariance matrices is also highly relevant in applications. Interpolation in this manifold is used in \cite{massart2019interpolation} for parametric order reduction and in \cite{bakker2018dynamic} for community detection. Optimization algorithms for various problems are proposed in \cite{bonnabel2010adaptive,meyer2011regression,marchand2016riemannian,mishra2011low}, including distance learning, distance matrix completion, and role-model extraction. Other studies have tackled the task of fitting curves to points on the manifold: \cite{gousenbourger2017piecewise} fit curves through covariance matrices that model a wind field; \cite{li2014conformational} use geodesic interpolation between Gram matrices to create intermediate protein conformations; \cite{kacem2018novel} employ Gram-matrix interpolation to track facial-expression changes in video data.

The set of $n$-dimensional covariance matrices is defined as
\begin{equation*}
	\Cov(n):=\left\{\Sigma \in \Sym(n) \,\vert\, \Sigma \succeq 0\right\},
\end{equation*}
where $\Sym(n)$ denotes the space of symmetric $n\times n$ matrices and the condition $\Sigma \succeq 0$ means that $\Sigma$ is positive semidefinite. The set $\Cov(n)$ forms a closed convex cone in the Euclidean space of square $n$-dimensional matrices $\Mat(n)$. It can be viewed as a stratified space whose strata are given by the smooth manifolds of covariance matrices of fixed rank $k$ (see \cite{thanwerdas2023bures}),
\begin{equation*}
	\Sym^+(n,k):=\left\{\Sigma \in \Cov(n) \,\vert\, \rk(\Sigma)=k\right\}.
\end{equation*}

The purpose of this work is to gain new insights on the Riemannian
manifold $(\Sym^+(n,k),d^{BW})$ of symmetric positive semidefinite matrices of fixed rank $k$, where $d^{BW}$ denotes the Bures--Wasserstein distance, through the diffeomorphism
\begin{equation}
	\label{eq:proposed model}
	\Sym^+(n,k)\;\cong\;
	M(n,k):=\St(n,k)\times_{O(k)}\Sym^{+}(k)
	=(\St(n,k)\times\Sym^{+}(k))/O(k).
\end{equation}
Here $\St(n,k)$ denotes the Stiefel manifold, $O(k)$ the orthogonal group and $M(n,k)$ the bundle associated
to the principal $O(k)$-bundle
$\pi_G\colon\St(n,k)\to\Gr(k,n)$. This associated-bundle model has been used in several applied works on low-rank covariance and subspace estimation (see, e.g., \cite{bonnabel2010adaptive, massart2019interpolation,bonnabel2010riemannian,smith2005covariance}), typically endowed with product-type Riemannian metrics. However, to the best of our knowledge, a systematic analysis of this model equipped with the Bures--Wasserstein metric---and its interplay with the Grassmannian geometry of $\Gr(k,n)$---has not been carried out.

From an information-geometric viewpoint, $\Cov(n)$ can be identified with the space of centred Gaussian measures on $\R^n$, with $\Sigma\in\Cov(n)$ corresponding to $\mathcal N(0,\Sigma)$. The Bures--Wasserstein distance $d^{BW}$ then coincides with the $L^2$-Wasserstein distance between these Gaussian measures \cite{malago2018wasserstein, villani2008optimal}, and provides an alternative to the classical Fisher--Rao metric on Gaussian families. The fixed-rank strata $\Sym^+(n,k)$ correspond to degenerate Gaussian measures supported on $k$-dimensional subspaces; understanding their Bures--Wasserstein geometry is therefore beneficial for structured covariance models such as factor models and low-rank approximations.

We recall that the Bures--Wasserstein distance is defined by
\begin{equation*}
	d^{BW}(\Sigma,\Lambda)^2=\Tr(\Sigma)+\Tr(\Lambda)-2\Tr\big((\Sigma^{1/2}\Lambda\Sigma^{1/2})^{1/2}\big).
\end{equation*}
With this distance, $\Cov(n)$ becomes a complete geodesic space (i.e. a complete metric space where the distance between two points is equal to the length of a minimizing geodesic connecting the two points) of nonnegative curvature (\cite{thanwerdas2023bures,takatsu2011wasserstein}).

The set of $n$-dimensional covariance matrices can be identified with the quotient space $\Mat(n) /O(n)$. The orbits of the quotient space are defined by the right action $\Mat(n) \times O(n) \ni (X,G) \mapsto XG \in \Mat(n)$ and the identification is given by the map 
\begin{equation}
	\label{int. identification map}
	\Mat(n) \ni X \mapsto X X^\top  \in \Cov(n).
\end{equation}

For a stratum $\Sym^+(n,k)$ a parallel construction to the one of $\Cov(n)$ can be built. In this case, we focus on the right action $\Mat(n,k)^\ast \times O(k) \ni (X,G) \mapsto XG \in \Sym^+(n,k)$, where $\Mat(n,k)^\ast$ denotes the manifold of full rank $n \times k$ matrices, and the consequent quotient space $\Mat(n,k)^\ast/O(k)$ which is diffeomorphic to $\Sym^+(n,k)$ via the map $\Mat(n,k)^\ast \ni X \mapsto X X^\top  \in \Sym^+(n,k)$. Thanks to this identification, the Euclidean metric descends to the Bures--Wasserstein Riemannian metric on $\Sym^+(n,k)$. For a more in-depth analysis we refer the reader to \cite{massart2020quotient}.

Before defining the metric, we recall the definition of the tangent space to $\Sym^+(n,k)$ at $\Sigma$
\begin{equation*}
	T_\Sigma \Sym^+(n,k) = \{V \in \Sym(n) \,|\,P_\Sigma^\perp V P_\Sigma^\perp=0\},
\end{equation*}
where $P_\Sigma^\perp$ denotes the orthogonal projectors onto the nullspace of $\Sigma$.
Now, following what is done in \cite{thanwerdas2023bures}, we set for a moment $\Sigma=Q D Q^\top $ with $D=\diag(d_1,...,d_k)$ positive definite and $Q \in \St(n,k)$ and call $\mathtt{S}_{\Sigma,V}=Q \cl S_D(Q^\top  V Q) Q^\top $, where $\cl{S}_A(B)=:C$ denotes the unique solution to the Sylvester equation $B=AC+CA$ for a symmetric positive definite matrix $A$ and a symmetric matrix $B$. For $V, \, W \in T_\Sigma\Sym^+(n,k)$, the Bures--Wasserstein metric is given by
\begin{align*}
	g_\Sigma^{BW(n,k)}(V,W)=&\Tr( \mathtt{S}_{\Sigma,V} \Sigma \mathtt{S}_{\Sigma,W}) + \Tr(P_\Sigma^\perp V \Sigma^\dagger W )\\=:&	g_\Sigma^{\ver}(V,W)+	g_\Sigma^{\hor}(V,W),
\end{align*}
where $\Sigma^\dagger$ denotes the Moore-Penrose inverse of $\Sigma$. Furthermore, the authors point out that $S_{\Sigma,W}$ and the metric are independent of the chosen decomposition.

With the above geometry constructions, the strata of $\Cov(n)$ have been studied in the literature for some years. The first results concerned the full rank principal stratum, where starting from the quotient geometry it was possible to derive equations for the exponential map (\cite{malago2018wasserstein}), unique logarithm (\cite{bhatia2019bures,massart2020quotient}), parallel transport (\cite{thanwerdas2023n}) and curvature (\cite{takatsu2010wasserstein,takatsu2011wasserstein}, \cite{thanwerdas2023n}). For the geometry of the stratum $\Sym^+(n,k)$ of fixed rank covariance matrices we refer to \cite{massart2020quotient}, \cite{massart2019curvature}, \cite{thanwerdas2023bures}; in particular the quotient geometry was studied in \cite{massart2020quotient}, a formula for the curvature tensor and the sectional curvature, as well as curvature bounds, were expressed --thanks to O'Neill's equations (\cite{o1966fundamental})-- in \cite{massart2019curvature}, and finally the exponential map, its definition domain, cut time and logarithms together with their multiplicity were found in \cite{thanwerdas2023bures}. Most recently, we gained some insights on the geometry (and in particular on the geodesics) of the stratified space as a whole (\cite{thanwerdas2023bures}).

In this paper, instead, we work with the associated-bundle model \eqref{eq:proposed model} and pull back the Bures--Wasserstein metric from $\Sym^+(n,k)$ to $M(n,k)$. This allows us to decouple the geometry into a horizontal component living over the Grassmannian and a vertical component on the symmetric positive definite (SPD) fiber $\Sym^+(k)$, and to derive explicit formulas for the metric and the corresponding geodesic equations.

\medskip

\noindent\textbf{Main contributions.}
The main contributions of this paper can be summarized as follows.
\begin{itemize}
	\item[(i)] We provide a detailed analysis of the associated-bundle model
	\begin{equation*}
		\Sym^+(n,k) \cong \St(n,k)\times_{O(k)}\Sym^{+}(k)
	\end{equation*}
	under the Bures--Wasserstein metric. In particular, we characterize the tangent spaces, the horizontal--vertical splitting, and the pullback metric on $M(n,k)$. This model gives a clear separation between changes and movements in the support subspace and the covariance within that subspace and it provides a standpoint to develop statistical procedures that are adaptive to this split. We stress again that the diffeomorphic model is not novel per se, but it has never been studied under the Bures--Wasserstein metric and a rigorous geometric perspective helps us to develop the further results named below.
	\item[(ii)] Within this model, we derive a system of ordinary differential equations (ODEs) for Bures--Wasserstein geodesics in bundle coordinates on $M(n,k)$. In particular, we identify a simple matrix--valued first integral which can be interpreted as the momentum associated with the $O(k)$--symmetry of the model. This conservation law decouples part of the dynamics and
	will be exploited in low-dimensional examples to recover known closed-form
	expressions for fixed-rank Bures--Wasserstein geodesics. This ODE formulation is well suited for numerical integration and algorithms on low-rank covariance manifolds.
	\item[(iii)] Utilizing the system of geodesic equations, we prove that the fibers $M_P(n,k)$, $P\in\Gr(k,n)$, are totally geodesic submanifolds. Equivalently, Bures--Wasserstein geodesics with vertical initial velocity remain in a fixed fiber, and geodesics joining two matrices with the same image stay in that image subspace. This gives a transparent geometric explanation of the algebraic uniqueness conditions for geodesics within a stratum observed in \cite{thanwerdas2023bures} and may prove useful in the development of algorithms aware of fiber preserving flows, e.g. in constrained optimization.
	\item[(iv)] We establish a one-to-one correspondence between the index set of Bures--Wasserstein logarithms on $(\Sym^+(n,k),d^{BW})$ and the index set of Grassmannian logarithms on $(\Gr(k,n),d_E)$, where $d_E$ denotes the Euclidean metric on the Grassmannian. This clarifies the precise relationship between multiplicities of minimizing geodesics in the Grassmann and fixed-rank covariance geometries.
\end{itemize}

From the information-geometric perspective, these results highlight the special role of the Grassmannian $\Gr(k,n)$ as a ``base space'' controlling the directions along which the support of a degenerate Gaussian distribution can move, while the SPD fiber $\Sym^+(k)$ governs the covariance structure within a fixed support subspace. The associated-bundle picture and the ODE system developed here provide a natural framework in which to study structured Gaussian models and low-rank covariance estimation under the Bures--Wasserstein metric.

\medskip

The article is organized as follows. In Section \ref{sec. preliminaries}, we introduce the needed geometrical framework. In particular, in Subsection \ref{subsec. fiber bundles} we recall basic definitions and useful results on fiber bundles, principal bundles, and associated bundles, while in Subsection \ref{subsec. principal bundle} we describe the principal fiber bundle $\pi_G: \St(n,k) \to \Gr(k,n)$.  

In Section \ref{sec. sym geometry}, we study the geometry of the stratum $\Sym^+(n,k)$. In Subsection \ref{subsec. diffeomoprhism}, we introduce and prove the diffeomorphism $\Sym^+(n,k) \cong M(n,k)$; in Subsection \ref{subsec. tangent spaces} we characterize the tangent spaces to $M(n,k)$; in Subsection \ref{subsec. metric on sym}, we study how the Bures--Wasserstein metric on $\Sym^+(n,k)$ translates into a metric on $M(n,k)$; in Subsection \ref{subsec. geodesics}, we analyze geodesics in $M(n,k)$, derive a system of geodesic equations and prove that the fibers are totally geodesic; and in Subsection \ref{subsec. bijection} we prove a bijection between the logarithms of $\Sym^+(n,k)$ and those of $\Gr(k,n)$.

Finally, for the sake of completeness and to preserve the flow of the exposition, proofs that are standard or particularly technical are collected in Appendix \ref{sec:appendix}.

\subsection*{Acknowledgments}
The author thanks professor Enrico Bernardi for valuable guidance and many helpful discussions during the preparation of this work.

	\section{Preliminaries}
	\label{sec. preliminaries}
	
	\subsection{Fiber bundles}
	\label{subsec. fiber bundles}
	
	Here we introduce the concepts from differential geometry that we will need in later Sections. We refer to \cite{sharpe2000differential,husemoller1966fibre,lee2000introduction} for a comprehensive treatment.
	
	\begin{definition}
		Let $E$, $B$ and $F$ be smooth manifolds and $\pi:E \to B$ a smooth surjection. The quadruple $(E,B,\pi,F)$ is called a \emph{smooth fiber bundle} with fiber $F$ if for every $b \in B$ there exists a neighbourhood $U \subset B$ and a diffeomorphism
		\begin{equation*}
		\psi_U: \pi^{-1}(U) \to U \times F
		\end{equation*}
		such that the following diagram commutes:
		\begin{equation*}
			\begin{tikzcd}
				\pi^{-1}(U) \arrow{r}{\psi_U} \arrow[swap]{dr}{\pi} & U \times F \arrow{d}{\pr_1} \\
				& U
			\end{tikzcd}
		\end{equation*}
		where $\pr_1$ denotes the canonical projection onto $U$. The manifold $B$ is called the \emph{base space}, $E$ the \emph{total space}, and $\pi$ the \emph{projection}. For $b \in B$, the set $E_b:=\pi^{-1}(b)$ is diffeomorphic to $F$ and is called the \emph{fiber} over $b$. The family $\{(U,\psi_U)\}$ is called a (smooth) \emph{local trivialization} of the bundle.
		
		Fiber bundles will often be denoted simply by $(E,\pi)$ or by $\pi: E \to B$.
	\end{definition}
	
	\begin{definition}
		Let $\pi_1: E_1 \to B$ and $\pi_2 : E_2 \to B$ be two fiber bundles over the same base space. A smooth map $\phi:E_1 \to E_2$ is called a \emph{bundle morphism} if it respects the projections, i.e.
		\begin{equation*}
			\pi_2 \circ \phi = \pi_1.
		\end{equation*}
	\end{definition}
	
	In particular, a bundle morphism $\phi: E_1 \to E_2$ maps each fiber $\pi_1^{-1}(b)$ into the corresponding fiber $\pi_2^{-1}(b)$.
	
	\begin{definition}
		A bundle morphism $\phi: E_1 \to E_2$ is called an \emph{isomorphism between fiber bundles} or a \emph{bundle isomorphism} if it is a diffeomorphism. The two bundles $\pi_1: E_1 \to B$ and $\pi_2 : E_2 \to B$ are then said to be \emph{isomorphic}.
	\end{definition}
	
	\begin{definition}
		Let $G$ be a Lie group. A \emph{principal $G$-bundle} is a fiber bundle $\pi:P \to B$ together with a smooth right action
		\begin{equation*}
		R_g: P \to P,\quad x \mapsto xg
		\end{equation*}
		of $G$ on $P$ such that:
		\begin{itemize}
			\item[(i)] for each $g \in G$, $R_g$ is a diffeomorphism;
			\item[(ii)] the action preserves the fibers and is free and transitive on them, i.e.\ for every fiber $P_b=\pi^{-1}(b)$ and every pair of points $x,\,y \in P_b$ there exists a unique $g \in G$ such that $y=xg$.
		\end{itemize}
	\end{definition}
	
	These properties imply that every fiber $P_b$ is diffeomorphic to $G$ and can be identified with $G$ itself.
	
	\begin{definition}
		\label{def:associated-bundle}
		Let $\pi:P \to B$ be a principal $G$-bundle. Let $F$ be a smooth manifold equipped with a left action of $G$ given by
		\begin{equation*}
		L_g: F \to F,\quad f \mapsto g \cdot f.
		\end{equation*}
		Define a right $G$-action on $P \times F$ by
		\begin{equation*}
			(p,f)\cdot g := (pg, g^{-1}\cdot f).
		\end{equation*}
		Denote by $E = P \times_G F := (P \times F)/G$ the corresponding orbit space, and write $[p,f]$ for the equivalence class of $(p,f)$. The projection
		\begin{equation*}
		\Tilde{\pi}: E \to B,\quad [p,f] \mapsto \pi(p),
		\end{equation*}
		turns $E$ into a fiber bundle with base $B$ and fiber $F$, called the \emph{fiber bundle associated} to the principal $G$-bundle $P$ or the \emph{associated $G$-bundle}.
	\end{definition}

	\begin{proposition}
		\label{prop. tangent spaces of associated bundle}
		Let $\pi:P \to B$ be a principal $G$-bundle and $E=P \times_G F$ an associated bundle with fiber $F$ and left action $L_g f=g \cdot f$. Then the tangent spaces of $E$ are given by
		\begin{equation*}
			T_{[p,f]} E = \big(T_pP \times T_f F\big)\big/\big\{(p \xi, -L_{\xi}^\ast f)\,\vert\, \xi \in \mathfrak{g}\big\},
		\end{equation*}
		where $\mathfrak{g}$ is the Lie algebra of $G$ and
		\begin{equation*}
		L_{\xi}^\ast := \partial_t \big(L_{\exp(t\xi)}(f)\big)\big|_{t=0}
		\end{equation*}
		is the infinitesimal action of $\xi$ on $F$.
	\end{proposition}
	
	\begin{proof}
		This follows from \cite[Theorems~1.1--1.3]{szczarba1964tangent}. A self-contained proof is given in the Appendix.
	\end{proof}
	With the same notation as in Definition \ref{def:associated-bundle} and Proposition \ref{prop. tangent spaces of associated bundle}, denoting by $\textrm{pr}_1:P \times F \to P$ the canonical projection and by $q:P \times F \to E$ the quotient map, we get
	\begin{equation*}
		\tilde \pi \circ q = \pi \circ \textrm{pr}_1.
	\end{equation*}
	Differentiating at $(p,f) \in  P \times F$ yields
	\begin{equation*}
		d \tilde \pi([p,f]) \circ dq(p,f)= d\pi(p)\circ d\textrm{pr}_1(p,f)
	\end{equation*}
	and thus, for $[x,y]\in T_{[p,f]} E$
	\begin{equation}
		\label{eq:differential of associated projection}
		\begin{aligned}
			d\tilde \pi([p,f])([x,y])=&d\tilde \pi([p,f])\circ dq(p,f)(x,y)\\=&d\pi(p)\circ d\textrm{pr}_1(p,f)(x,y)\\=&d\pi(p)(x).
		\end{aligned}		
	\end{equation}
	This gives the vertical spaces
	\begin{equation}
		\label{eq:vertical spaces}		
			\cl{V}_{[p,f]} := \Ker\left(d\tilde \pi([p,f])\right)
			= \{[x,y] \in T_{[p,f]}E \,\vert\, x \in \Ker (d\pi(p))\}.
	\end{equation}

	We will also need the following Lemma, which allows us to upgrade fiberwise diffeomorphisms to a global bundle isomorphism.
	
	\begin{lemma}
		\label{lemma:local to global diff}
		Let $\pi_M: M \to B$ and $\pi_N: N \to B$ be two smooth fiber bundles over the same base space. Let $\phi: M \to N$ be a smooth bundle morphism 
		and assume that, for every $b \in B$, the restriction
		\begin{equation*}
		\phi_b := \phi\big|_{\pi_M^{-1}(b)} : \pi_M^{-1}(b) \to \pi_N^{-1}(b)
		\end{equation*}
		is a diffeomorphism. Then $\phi$ is a global diffeomorphism.
	\end{lemma}
	
	\begin{proof}
		This is a modification of \cite[Theorem~3.2]{husemoller1966fibre}. A complete proof is given in the Appendix.
	\end{proof}

	\subsection{The Grassmannian as a principal fiber bundle}
	\label{subsec. principal bundle}
	
	We now briefly review the well-known principal $O(k)$-bundle
	\begin{equation*}
		\pi_G:\St(n,k) \to \Gr(k,n), \qquad Q \mapsto \Span(Q),
	\end{equation*}
	where $\St(n,k)$ denotes the Stiefel manifold of $n\times k$ matrices with orthonormal columns, and $\Gr(k,n)$ the Grassmann manifold of $k$-dimensional subspaces of $\R^n$. For a comprehensive review we refer the reader to \cite{bendokat2024grassmann}.
	
	Let $Q \in \St(n,k)$. The tangent space to $\St(n,k)$ at $Q$ can be characterized as
	\begin{equation*}
		T_Q \St(n,k)=\{V \in \Mat(n,k) \,\vert\, Q^\top V + V^\top Q=0\}.
	\end{equation*}
	If we fix a matrix $Q_\perp \in \Mat(n,n-k)$ that completes $Q$ to an orthonormal basis of $\R^n$, then every $V \in T_Q \St(n,k)$ can be uniquely decomposed as
	\begin{equation*}
		V = Q A + Q_\perp B, \qquad A \in \mathfrak{o}(k)=\Skew(k),\; B \in \Mat(n-k,k),
	\end{equation*}
	where $\mathfrak{o}(k)$ denotes the Lie algebra of the orthogonal group $O(k)$ (see \cite[Section~2.2]{bendokat2024grassmann}).
	
	We identify the tangent space to $\Gr(k,n)$ at a point $P$ with
	\begin{equation*}
		T_P \Gr(k,n) \cong \Hom(P, P^\perp).
	\end{equation*}
	Once a basis is specified, we may further identify $\Hom(P, P^\perp) \cong \Mat(n-k,k)$.
	
	The vertical space at $Q \in \St(n,k)$ is given by
	\begin{equation*}
		\cl{V}_Q:= \Ker\big(d\pi_G(Q)\big)=\{QA \,\vert\, A \in \mathfrak{o}(k)\}.
	\end{equation*}
	
	The standard construction (see \cite{bendokat2024grassmann}) is to endow both $\St(n,k)$ and $\Gr(k,n)$ with the Euclidean (Frobenius) metric, i.e.
	\begin{alignat}{2}
		\label{metric on grassmannian}
		g^{\St}_Q(V_1,V_2)&=\Tr(V_1 V_2^\top), &\qquad& V_1 , V_2 \in T_Q \St(n,k), \nonumber \\
		g^{\Gr}_P(B_1,B_2)&=\Tr(B_1 B_2^\top), &\qquad& B_1, B_2 \in \Mat(n-k,k) \cong T_P \Gr(k,n),
	\end{alignat}
	where the Grassmannian metric is independent of the chosen matrix representatives of $P$.

	The \emph{horizontal} space at $Q$ is then defined as the orthogonal complement of $\cl{V}_Q$ in $T_Q\St(n,k)$ with respect to the Euclidean metric, namely
	\begin{equation*}
		\cl{H}_Q:= \cl{V}_Q^\perp=\{Q_\perp B \,\vert\, B \in \Mat(n-k,k)\}.
	\end{equation*}

	\begin{proposition}
		\label{prop:principal bundlle submersion}
		The projection $\pi_G:(\St(n,k),g^{\St}) \to (\Gr(k,n),g^{\Gr})$ is a Riemannian submersion, i.e.\ the restriction
		\[
		d\pi_G(Q)\big|_{\cl{H}_Q} : \cl{H}_Q \to T_P \Gr(k,n), \qquad P=\Span(Q),
		\]
		is an isometry for every $Q\in\St(n,k)$.
	\end{proposition}
	
	\begin{proof}
		This is shown in \cite{bendokat2024grassmann}. For the reader's convenience we recall the main argument and provide a proof in the Appendix. In particular, in that proof we derive an explicit formula for the differential of the projection:
		\begin{equation}
			\label{eq:diff of projection in principal bundle}
			d\pi_G(Q)(V)=B \in \Mat(n-k, k)
		\end{equation}
		when $V=Q_\perp B$ and $[Q \,\, Q_\perp]$ is used as an orthonormal basis of $\R^n$.
	\end{proof}

\section{The Bures--Wasserstein geometry of the fixed-rank stratum}
\label{sec. sym geometry}

\subsection{The diffeomorphism}
\label{subsec. diffeomoprhism}

As recalled in Subsection~\ref{subsec. principal bundle}, $\St(n,k)$ carries a principal $O(k)$-bundle structure with base space $\Gr(k,n)$ and structure group $O(k)$ acting by right multiplication. 

We equip $\Sym^+(k)$ with the smooth left $O(k)$-action
\begin{equation*}
G \cdot D := G D G^\top
\end{equation*}
and define a right $O(k)$-action on $\St(n,k) \times \Sym^+(k)$ by
\begin{equation}
	\label{eq:right-action}
	(Q,D)\cdot G:=(QG, G^\top D G).
\end{equation}
The action defined in \eqref{eq:right-action} is smooth and free; indeed, if $(Q,D) \cdot G=(Q,D)$ then $QG=Q$ and thus
\begin{equation*}
	G=Q^\top QG=Q^\top Q=I.
\end{equation*}
Moreover, in view of \cite[Corollary 21.6]{lee2000introduction}, compactness of $O(k)$ gives properness of the right action in \eqref{eq:right-action}.
We then form the associated fiber bundle
\begin{equation*}
	M(n,k)=\St(n,k) \times_{O(k)}  \Sym^+(k)
	=\big(\St(n,k)\times \Sym^+(k)\big)/O(k),
\end{equation*}
where the equivalence relation is
\begin{equation*}
	(Q,D)\sim (Q,D)\cdot G, \qquad G\in O(k),
\end{equation*}
and equivalence classes are denoted by $[Q,D]$. We define the quotient map
\begin{equation*}
	q_M:\St(n,k)\times \Sym^+(k) \to M(n,k), \quad q_M(Q,D)=[Q,D].	
\end{equation*} 
Applying \cite[Theorem 21.10]{lee2000introduction}, we find that, since the right action of $O(k)$ on $\St(n,k) \times \Sym^+(k)$ defined in \eqref{eq:right-action} is smooth free and proper, the quotient map $q_M$ is a smooth surjective submersion.
The projection of the associated bundle
\begin{equation*}
\pi_M: M(n,k) \to \Gr(k,n)
\end{equation*}
is defined by
\begin{equation*}
	\pi_M([Q,D])=\pi_G(Q)=\Span(Q),
\end{equation*}
and the fiber over $P\in\Gr(k,n)$ is
\begin{equation*}
	M_P(n,k)=\{[Q,D] \,\vert\, \Span(Q)=P\}.
\end{equation*}

It is a standard fact in the theory of associated fiber bundles that every fiber $M_P(n,k)$ is diffeomorphic to the model fiber $\Sym^+(k)$ (see, for example, \cite[Chapter~1]{kobayashi1963foundations} and \cite[Chapter~4]{husemoller1966fibre}).

On the other hand, we consider the bundle
\begin{equation*}
	\pi_{\Sym}: \Sym^+(n,k) \to \Gr(k,n), \qquad \pi_{\Sym}(\Sigma)=\Imm(\Sigma),
\end{equation*}
where $\Imm(\Sigma)$ denotes the image (column space) of $\Sigma$. The fiber over $P\in\Gr(k,n)$ is
\begin{equation*}
	\Sym_P^+(n,k)=\{\Sigma \in \Sym(n) \,\vert\, \Sigma \succeq 0,\ \rk(\Sigma)=k,\ \Imm(\Sigma)=P\}.
\end{equation*}
Thus $\Sym^+(n,k)$ is also a fiber bundle over $\Gr(k,n)$ with the same base space as $M(n,k)$.

We define the smooth map
\begin{equation*}
	\Tilde \phi: \St(n,k) \times \Sym^+(k) \to \Sym^+(n,k), \qquad \Tilde \phi(Q,D)=Q D Q^\top .
\end{equation*}
We observe that $\Tilde \phi(Q G, G^\top  D G)=\Tilde \phi(Q,D)$ for all $G \in O(k)$. Therefore, thanks to \cite[Theorem 4.30]{lee2000introduction}, $\Tilde \phi$ descends to a unique smooth map $\phi$
\begin{equation*}
	\phi: M(n,k) \to \Sym^+(n,k), \qquad \phi([Q,D])=\Tilde \phi(Q,D)=Q D Q^\top .
\end{equation*}
Moreover, $\phi$ preserves fibers:
\begin{equation*}
\pi_{\Sym}\big(\phi([Q,D])\big)
=\Imm(QDQ^\top )
=\Span(Q)
=\pi_M([Q,D]).
\end{equation*}
Hence $\phi$ is a smooth bundle morphism from $M(n,k)$ to $\Sym^+(n,k)$ over $\Gr(k,n)$.

\begin{proposition}
	\label{prop:fiberwise-diffeo}
	The smooth bundle morphism $\phi: M(n,k) \to \Sym^+(n,k)$ is a fiberwise diffeomorphism, i.e., the restriction
	\begin{equation*}
		\phi|_{M_p(n,k)}: M_p(n,k) \to \Sym_P^+(n,k)
	\end{equation*}
	is a diffeomorphism.
\end{proposition}

\begin{proof}
	Fix $P\in\Gr(k,n)$ and choose $Q_0\in\St(n,k)$ such that $\Span(Q_0)=P$.
	If $[Q,D]\in M_P(n,k)$, then $Q$ and $Q_0$ span the same subspace, hence there exists
	$G\in O(k)$ such that $Q=Q_0G$. Therefore
	\[
	[Q,D]=[Q_0G,D]=[Q_0,G^\top D G],
	\]
	so every class in $M_P(n,k)$ admits a representative of the form $[Q_0,D]$ with $D\in\Sym^+(k)$.
	This yields the identification $M_P(n,k)\simeq \Sym^+(k)$ via $[Q_0,D]\leftrightarrow D$.
	
	Consider the restriction of $\phi$ to the fiber and, under the above identification, define
	\[
	\phi_{M_P(n,k)}:\Sym^+(k)\to \Sym_P^+(n,k),\qquad \phi_{M_P(n,k)}(D)=Q_0 D Q_0^\top.
	\]
	The map $D\mapsto Q_0 D Q_0^\top$ is the restriction of a linear map
	$\Mat(k)\to\Mat(n)$, hence it is smooth. It is injective because $Q_0$ has full column rank.
	It is surjective onto $\Sym_P^+(n,k)$ by definition of $\Sym_P^+(n,k)$: if $\Sigma\in\Sym_P^+(n,k)$,
	then $\Imm(\Sigma)=P=\Span(Q_0)$ and $\Sigma$ is symmetric positive semidefinite of rank $k$, hence
	$\Sigma$ is positive definite on $P$ and admits the representation $\Sigma=Q_0 D Q_0^\top$ with
	\[
	D:=Q_0^\top \Sigma Q_0\in\Sym^+(k).
	\]
	Indeed, $D$ is symmetric, and for any $x\neq 0$ one has
	\[
	x^\top D x = (Q_0x)^\top \Sigma (Q_0x)>0,
	\]
	since $Q_0x\in P\setminus\{0\}$.
	
	Finally, the inverse map is explicitly
	\[
	\phi_{M_P(n,k)}^{-1}:\Sym_P^+(n,k)\to \Sym^+(k),\qquad \phi_{M_P(n,k)}^{-1}(\Sigma)=Q_0^\top \Sigma Q_0,
	\]
	which is again the restriction of a linear map, hence smooth. Therefore $\phi_{M_P(n,k)}$ is a
	diffeomorphism.
\end{proof}
We can now apply Lemma~\ref{lemma:local to global diff} to the bundle morphism $\phi: M(n,k) \to \Sym^+(n,k)$.

\begin{theorem}
	\label{thm:global-diffeo}
	Let
	\begin{equation*}
		M(n,k)=\St(n,k) \times_{O(k)}  \Sym^+(k)
		=\big(\St(n,k)\times \Sym^+(k)\big)/O(k).
	\end{equation*}
	Then the map
	\begin{equation*}
	\phi: M(n,k) \xrightarrow{\ \cong\ } \Sym^+(n,k),\qquad [Q,D]\mapsto Q D Q^\top ,
	\end{equation*}
	is a global diffeomorphism of fiber bundles over $\Gr(k,n)$. In particular, each fiber of $\Sym^+(n,k)$ is diffeomorphic to $\Sym^+(k)$.
\end{theorem}

Practically, we can now assign in a smooth manner to a matrix $\Sigma=Q D Q^\top $ (with $Q\in\St(n,k)$ and $D\in\Sym^+(k)$) a class $[Q,D]\in M(n,k)$, and conversely, any class $[Q,D]$ corresponds to the fixed-rank covariance matrix $\Sigma=Q D Q^\top $.

\subsection{Tangent spaces}
\label{subsec. tangent spaces}

We now identify the tangent spaces of $M(n,k)$ by means of Proposition~\ref{prop. tangent spaces of associated bundle}. In our case, the left action of $O(k)$ on $\Sym^+(k)$ is given by
\begin{equation*}
L_G(D)=G D G^\top .
\end{equation*}
For $A \in \mathfrak{o}(k)$ and $D\in\Sym^+(k)$, the infinitesimal action is
\begin{align*}
	L_A^\ast(D)
	&= \partial_t\big(\exp(tA) \, D\, \exp(tA)^\top\big)\big|_{t=0} \\
	&= AD -  DA.
\end{align*}
Hence, by Proposition~\ref{prop. tangent spaces of associated bundle}, the subspace that is quotiented out in $T_{[Q,D]}M(n,k)$ is
\begin{equation*}
\{(Q A,-L_A^\ast(D)) \,\vert\, A \in \mathfrak{o}(k)\},
\end{equation*}
and the equivalence relation on $T_Q\St(n,k)\times T_D\Sym^+(k)$ is
\begin{equation}
	\label{eq:equivalence relation}
(V,T)\sim (V+Q A,T-L_A^\ast(D)), \qquad A\in\mathfrak{o}(k).
\end{equation}

Recall that every $V \in T_Q \St(n,k)$ can be uniquely decomposed as
\begin{equation*}
V=Q A + Q_\perp B, \qquad A \in \mathfrak{o}(k),\; B \in \Mat(n-k,k),
\end{equation*}
for a suitable completion $Q_\perp$ of $Q$ to an orthonormal basis of $\R^n$. Using the equivalence relation in \eqref{eq:equivalence relation}, any pair $(V,T)$ is equivalent to a unique representative of the form $(Q_\perp B,\widetilde{T})$, obtained by absorbing the $QA$ component into the vertical part via $T\mapsto T-L_A^\ast(D)$.

Thus we obtain the following description:

\begin{align*}
	T_{[Q,D]}M(n,k)
	= &\{[Q_\perp B, T] \,\vert\, B \in \Mat(n-k,k),\ T \in \Sym(k)\}
	\\\cong &\{(Q_\perp B, T)\,\vert\, B \in \Mat(n-k,k),\ T \in \Sym(k)\},
\end{align*}
where, for notational convenience, we identify each class with its representative $(Q_\perp B, T)$.

Now, Equation (\ref{eq:differential of associated projection}) can be applied to this case, yielding
\begin{align*}
	d\pi_M([Q,D])([Q_\perp B, T])= d\pi_G(Q)(Q_\perp B)= B.
\end{align*}

This yields the horizontal--vertical splitting
\begin{align*}
	\cl{V}_{[Q,D]} &= \Ker\big(d\pi_M([Q,D])\big)
	= \{[0, T] \,\vert\, T \in \Sym(k)\},
	 \cong\{(0, T) \,\vert\, T \in \Sym(k)\},\\
	\cl{H}_{[Q,D]}&= \{[Q_\perp B, 0] \,\vert\, B \in \Mat(n-k,k)\}\cong\{(Q_\perp B, 0) \,\vert\, B \in \Mat(n-k,k)\}.
\end{align*}
In other words, in the bundle coordinates $[Q,D]$ with tangent representatives $(Q_\perp B,T)$, horizontal directions are encoded by $B$ (moving the subspace in the Grassmannian) and vertical directions by $T$ (changing the covariance within a fixed subspace).

\subsection{Pullback metric}
\label{subsec. metric on sym}

Let $(N,g)$ be a Riemannian manifold and let $f \colon M \to N$ be a smooth map. 
The \emph{pullback metric} $f^\ast g$ on $M$ is defined by
\begin{equation*}
(f^\ast g)_x(V,W) := g_{f(x)}\big(df_x(V), df_x(W)\big),
\qquad x\in M,\; V,W\in T_xM.
\end{equation*}
In our setting we take $N = \Sym^+(n,k)$ endowed with the Bures--Wasserstein metric $g^{BW(n,k)}$, 
and $ M(n,k)$ equipped with the diffeomorphism
\begin{equation*}
\phi \colon M(n,k) \xrightarrow{\ \cong\ } \Sym^+(n,k),\qquad [Q,D] \mapsto Q D Q^\top 
\end{equation*}
from Theorem~\ref{thm:global-diffeo}. 
We then define the Riemannian metric $h$ on $M(n,k)$ as the pullback
\begin{equation*}
h := \phi^\ast g^{BW(n,k)}.
\end{equation*}

\begin{proposition}
	\label{prop:pullback-metric}
	Let $[Q,D]\in M(n,k)$ with $Q\in\St(n,k)$, $D\in\Sym^+(k)$, and let 
	$Q_\perp\in\St(n,n-k)$ be such that $[Q \, Q_\perp] \in O(n)$. Moreover, let $V, \, W \in T_{[Q,D]}M(n,k)$ be written as
	\begin{equation*}
	V := [Q_\perp B_V, T_V],\qquad
	W := [Q_\perp B_W, T_W],
	\end{equation*}
	with $B_V,B_W\in\Mat(n-k,k)$ and $T_V,T_W\in\Sym(k)$.
	Then the pullback metric $h = \phi^\ast g^{BW(n,k)}$ on $M(n,k)$ satisfies
		\begin{equation}\label{eq:pullback-metric}
			h_{[Q,D]}(V,W)
			= \Tr\big(B_V D B_W^\top \big)
			+ \Tr\big(\cl S_D(T_V)\,D\,\cl S_D(T_W)\big).
		\end{equation}
	In particular, $h$ is a well-defined $O(k)$-invariant Riemannian metric on $M(n,k)$, and $\phi \colon (M(n,k),h)\to (\Sym^+(n,k),g^{BW(n,k)})$ is a Riemannian isometry.
\end{proposition}
\begin{proof}
	The proof is shown in the Appendix.
\end{proof}

\begin{remark}
	We observe that the pullback metric naturally splits into an  horizontal and a vertical part
	\begin{gather*}
		h_{[Q,D]}^{\hor}(V,W):=\Tr\big(B_V D B_W^\top \big),\\
		h_{[Q,D]}^{\ver}(V,W):=\Tr\big(\cl S_D(T_V)\,D\,\cl S_D(T_W)\big),
	\end{gather*}
	and that the horizontal metric retains vertical information, i.e., it is influenced  by vertical position coordinate $D$. 
	
	Moreover, the Moore-Penrose inverse contribution $\Sigma^\dagger$ in the horizontal part of the metric vanishes. The detailed computation is carried out in the proof of Proposition~\ref{prop:pullback-metric}.
	We briefly recall here the mechanism behind this simplification.
	Writing $\Sigma=QDQ^\top$ gives $P^\perp_\Sigma=Q_\perp Q_\perp^\top$ and $\Sigma^\dagger=QD^{-1}Q^\top$.
	For $V=[Q_\perp B_V,T_V]$ we have $d\phi(V)=Q_\perp B_V DQ^\top+QDB_V^\top Q_\perp^\top+QT_VQ^\top$, hence
	$P^\perp_\Sigma\,d\phi(V)=Q_\perp B_VDQ^\top$ and
	$P^\perp_\Sigma\,d\phi(V)\,\Sigma^\dagger=Q_\perp B_VQ^\top$ (the factor $D$ cancels with $D^{-1}$).
	Similar computations are performed for $d\phi(W)$. Consequently $g^{\mathrm{hor}}_\Sigma(d\phi(V),d\phi(W))=\Tr(B_VDB_W^\top)$, where the terms
	involving $B_W$ and $T_W$ vanish by $Q^\top Q_\perp=0$ and by cyclicity of the trace.
\end{remark}

\subsection{Geodesics}
\label{subsec. geodesics}
In this subsection we study geodesics on the associated-bundle model
\[
M(n,k)=\St(n,k)\times_{O(k)}\Sym_+(k)
\]
endowed with the pullback Bures--Wasserstein metric
\[
h=\phi^*g_{BW(n,k)}.
\]
Since $\phi:(M(n,k),h)\to (\Sym_+(n,k),g_{BW(n,k)})$ is a Riemannian isometry, any
description of geodesics in $M(n,k)$ immediately yields the corresponding description
on the fixed-rank covariance manifold $\Sym_+(n,k)$.

We first derive a geodesic system of differential equations in bundle coordinates and recover a first integral. As a corollary, we show that the fibers
are totally geodesic. We then turn to minimizing logarithms and prove that the
multiplicity of minimizing Bures--Wasserstein logarithms on $\Sym_+(n,k)$ is governed
by the same orthogonal degree of freedom as for minimizing logarithms on the
Grassmannian $\Gr(k,n)$.

\subsubsection{ODE system for geodesic equations and totally geodesic fibers}
\label{subsec:ode_system}

We now derive the geodesic equations on $M(n,k)$ via a variational approach. 

\medskip

\begin{theorem}[Geodesic ODE system]
	\label{thm:BW-geodesic-ODE}
	Let $M(n,k)=\St(n,k)\times_{O(k)}\Sym^+(k)$ be endowed with the pullback Bures--Wasserstein metric $h=\phi^\ast g^{BW(n,k)}$. Consider a smooth curve
	\begin{equation*}
	\gamma(t) = [Q(t),D(t)] \in M(n,k), \qquad t\in I\subset\R,
	\end{equation*}
	with $Q(t)\in\St(n,k)$ and $D(t)\in\Sym^+(k)$. Let $Q_\perp(t)$ be any smooth choice of orthonormal complement of $Q(t)$ so that
	\begin{equation*}
	\mathbf{Q}(t) := \big[Q(t)\;\;Q_\perp(t)\big] \in O(n),
	\end{equation*}
	and write the velocity as
	\[
	\dot\gamma(t) = (Q_\perp(t) B(t), T(t)),
	\]
	with $B(t)\in\Mat(n-k,k)$ and $T(t)\in\Sym(k)$. Let $S(t):= \cl{S}_{D(t)}\left(T(t)\right)$.
	Then $\gamma$ is a Bures--Wasserstein geodesic in $M(n,k)$ if and only if it satisfies the bundle coordinate system
	\begin{equation}
		\label{eq:full-geodesic-system}
		\begin{cases}
			\dot{\mathbf{Q}} = \mathbf{Q}\,\mathbf{B}, & \mathbf{Q}(0)=\mathbf{Q}_0,\\[0.25em]
			\dot D = D S + S D, & D(0)=D_0,\\[0.25em]
			\dot B = -B(DS+SD)D^{-1}, & B(0)=B_0,\\[0.25em]
			\dot S = B^\top  B - S^2, & S(0)=S_{D_0}(\dot D(0))=S_0,
		\end{cases}
	\end{equation}
	where
	\[
	\mathbf{B}(t) :=
	\begin{bmatrix}
		0 & -B(t)^\top \\
		B(t) & 0
	\end{bmatrix} \in \mathfrak{o}(n)=\Skew(n).
	\]
	In particular, the first line of \eqref{eq:full-geodesic-system} is equivalent to
	\[
	\dot Q = Q_\perp B,
	\qquad
	\dot Q_\perp = -Q B^\top ,
	\]
	so that the columns of $Q(t)$ and $Q_\perp(t)$ remain orthonormal for all $t$.
\end{theorem}

\begin{proof}
	The proof is shown in the Appendix.
\end{proof}

\medskip

\begin{remark}[Conserved momentum]
	\label{rem:conserved-momentum}
	As highlighted in the proof of Theorem~\ref{thm:BW-geodesic-ODE} in the Appendix,
	the Lagrangian
	\[
	L(B,D,T) = \frac12 \Tr(BDB^\top)
	+ \frac12 \Tr\bigl(\cl S_D(T)\, D\, \cl S_D(T)\bigr)
	\]
	does not depend on the base point $Q \in \St(n,k)$. In particular, the
	conjugate momentum with respect to the horizontal velocity $B$,
	\[
	\partial_B L = BD \in \Mat(n-k,k),
	\]
	is conserved along Bures--Wasserstein geodesics, i.e.
	\begin{equation}
		\label{eq:conserved-momentum}
	\frac{\mathrm d}{\mathrm dt}\bigl(B(t)D(t)\bigr) = 0.
	\end{equation}
	Equivalently, the matrix $K := BD$ is a first integral of the geodesic flow
	on $M(n,k)$, and can be viewed as the momentum associated with the $O(k)$-symmetry
	of the bundle model.
	This feature suggests the design of structure-preserving numerical schemes
	that exactly maintain $BD$, in analogy with angular-momentum-preserving
	integrators in classical mechanics. In the rank-one examples below this conservation law
	reduces the geodesic system to a one-dimensional problem and allows us to
	recover the known fixed-rank Bures--Wasserstein geodesics in closed form.
\end{remark}

We can utilize Theorem \ref{thm:BW-geodesic-ODE} to show that the fibers $M_P(n,k), \, P \in \Gr(k,n)$ are totally geodesic.

	\begin{corollary}[Totally geodesic fibers]
	\label{prop:fibers-totally-geodesic}
	Fix $P\in \Gr(k,n)$. The fiber $M_P(n,k)=\pi_M^{-1}(P)$ is totally geodesic in $(M(n,k),h)$.
\end{corollary}

\begin{proof}
	Let $\gamma(t)=[Q(t),D(t)],\, t\in I \subset \R$ be a geodesic with $\gamma(0)\in M_P(n,k)$ and $\dot\gamma(0)$ vertical.
	In bundle coordinates $\dot\gamma(t)=(Q_\perp(t)B(t),T(t))$, verticality at $t=0$ means $B(0)=0$.
	By Theorem~\ref{thm:BW-geodesic-ODE}, $(Q,D,B,S)$ satisfies the geodesic ODE system, and in particular
	\[
	\dot B(t)=-B(t)(D(t)S(t)+S(t)D(t))D(t)^{-1}.
	\]
	Thus $B\equiv 0$ is a solution with $B(0)=0$. By uniqueness for smooth ODEs, $B(t)\equiv 0$ for all $t \in I$.
	Consequently $\dot Q(t)=Q_\perp(t)B(t)=0$, so $Q(t)\equiv Q(0)$ and $\pi_M(\gamma(t))=\Span(Q(t))$ is constant.
	Hence $\gamma(t)\in M_P(n,k)$ for all $t$, proving that $M_P(n,k)$ is totally geodesic.
\end{proof}

\begin{remark}
	Corollary \ref{prop:fibers-totally-geodesic} is consistent with the fixed-rank Bures--Wasserstein geometry described in \cite{thanwerdas2023bures}, and in fact provides a more geometric perspective on it. In particular, \cite[Theorem~5.2]{thanwerdas2023bures} states that the Bures--Wasserstein geodesic connecting two matrices $\Sigma,\Lambda\in\Sym^+(n,k)$ is unique if and only if $\rk(\Sigma\Lambda)=\rk(\Sigma)$, i.e.\ if and only if $\Sigma$ and $\Lambda$ have the same image subspace. In our bundle picture this means that the two points lie in the same fiber $M_P(n,k)$, and Corollary~\ref{prop:fibers-totally-geodesic} shows that the geodesic joining them remains inside this fiber. Thus the uniqueness criterion of \cite{thanwerdas2023bures} acquires a natural geometric interpretation, and, in addition, such instances can be studied within the full-rank Bures--Wasserstein framework on the fiber, where the geometry is technically and computationally more tractable.
\end{remark}

We now specialize the ODE system of Theorem \ref{thm:BW-geodesic-ODE} to low-dimensional examples where we can solve it in closed form and explicitly recover the known Bures--Wasserstein geodesics on $\Sym^+(n,k)$.

\begin{example}
	\label{ex. n=2 k=1}
	Let $n=2$ and $k=1$. Then $\St(2,1)\cong S^1$ and $O(1)=\{\pm1\}$, so
	\begin{equation*}
	\Gr(1,2)=\R \P^1 \cong S^1/\{\pm 1\},
	\end{equation*}
	where $\R \P^1$ denotes the real projective line. Moreover $\Sym^+(1)=(0,+\infty)$, and therefore
	\begin{equation*}
	M(2,1) = (S^1 \times (0,+\infty))/\{\pm 1\}, \qquad
	\pi_M: M(2,1) \to \R P^1,
	\end{equation*}
	with the equivalence relation $(q,d)\sim(-q,d)$. The global diffeomorphism of Theorem~\ref{thm:global-diffeo} is
	\begin{equation*}
	\phi: M(2,1) \xrightarrow{\ \cong\ } \Sym^+(2,1), \qquad [q,d] \mapsto d\,q q^\top .
	\end{equation*}
	
	We parametrize $S^1$ by
	\begin{equation*}
	q(\theta)=\begin{pmatrix}
		\cos\theta\\
		\sin\theta
	\end{pmatrix},
	\qquad \theta \in(-\pi,\pi],
	\end{equation*}
	and choose the perpendicular vector
	\begin{equation*}
	q_\perp(\theta)=\begin{pmatrix}
		-\sin\theta\\
		\cos\theta
	\end{pmatrix}.
	\end{equation*}
	In these coordinates, a point in $M(2,1)$ is represented by $[\theta,d]$ with the identification $(\theta,d)\sim(\theta+\pi,d)$, and a tangent vector at $[\theta,d]$ is represented by $(q_\perp(\theta)b,u)$ with $b,u\in\R$.
	
	In one dimension the Sylvester equation reduces to
	\begin{equation*}
	d s + s d = u \quad\Longrightarrow\quad s = \frac{u}{2d}.
	\end{equation*}
	The metric from Subsection~\ref{subsec. metric on sym} becomes
	\begin{equation*}
	h_{[\theta,d]}(q_\perp b, u)= d\,b^2 + d\,s^2 = d\,b^2 + \frac{u^2}{4d}.
	\end{equation*}
	The geodesic system \eqref{eq:full-geodesic-system} reduces to
	\begin{equation}
		\label{ex. geodesic system}
		\begin{cases}
			\dot\theta = b, & \theta(0)=\theta_0,\\[0.2em]
			\dot d = 2ds, & d(0)=d_0,\\[0.2em]
			\dot b = -2bs, & b(0)=b_0,\\[0.2em]
			\dot s = b^2 - s^2, & s(0)=s_0.
		\end{cases}
	\end{equation}
	From \eqref{eq:conserved-momentum} we have $bd=b_0 d_0$ along the geodesic.
	
	From \eqref{ex. geodesic system} we immediately see that the fibers are totally geodesic: if $b_0=0$, then $b\equiv 0$ and hence $\dot\theta\equiv0$, so $\theta$ is constant and the geodesic remains in the fiber over a fixed point of $\R P^1$. In that case, the system reduces to
	\begin{equation*}
	\dot d = 2ds,\qquad \dot s = -s^2,
	\end{equation*}
	whose solution is
	\begin{equation*}
	d(t) = d_0(1+s_0 t)^2.
	\end{equation*}
	This matches the one-dimensional Bures--Wasserstein geodesic on $\Sym^+(1)$ obtained as a special case of \cite[Theorem~4.2]{thanwerdas2023bures}.
	
	For $b_0\neq0$ we can still solve \eqref{ex. geodesic system} explicitly. Define $r=s/b$ and compute
	\begin{equation*}
	\dot r = \frac{\dot s\,b - \dot b\,s}{b^2} = b(r^2+1),
	\qquad
	\dot b = -2b^2 r.
	\end{equation*}
	Dividing these equations yields
	\begin{equation*}
	\frac{db}{dr} = \frac{\dot b}{\dot r} = -\frac{2rb}{r^2+1},
	\end{equation*}
	or equivalently
	\begin{equation*}
	\frac{d (\ln b)}{dr} = -\frac{2r}{r^2+1},
	\end{equation*}
	which integrates to
	\begin{equation}
		\label{ex. b r relationship}
		b(r^2+1)=b_0(r_0^2+1) =: C,
	\end{equation}
	where $r_0=s_0/b_0$. Integrating $\dot r = b(r^2+1)$ then gives
	\begin{equation*}
	r(t) = r_0 + Ct.
	\end{equation*}
	Using \eqref{ex. b r relationship} we obtain
	\begin{equation*}
	b(t) = \frac{C}{(Ct+r_0)^2+1},
	\qquad
	s(t) = r(t)b(t) = \frac{C(Ct+r_0)}{(Ct+r_0)^2+1}.
	\end{equation*}
	From $bd=b_0d_0$ we also have
	\begin{equation}
		\label{ex. d equation}
		d(t)=d_0\,\frac{(Ct+r_0)^2+1}{r_0^2+1}.
	\end{equation}
	
	Finally, we integrate $\dot\theta=b$:
	\begin{equation*}
	\theta(t)=\theta_0 + \int_0^t b(\tau)\,d\tau
	=\theta_0 + \big(\arctan(Ct+r_0)-\arctan(r_0)\big).
	\end{equation*}

	To compare with the known expressions for Bures--Wasserstein geodesics on $\Sym^+(2,1)$, one writes $q(\theta(t))$ in the orthonormal basis $\{q_0,w_0\}:=\{q(\theta_0),q_\perp(\theta_0)\}$, and using elementary trigonometric identities together with \eqref{ex. d equation}, checks that
	\begin{equation*}
	\sqrt{d(t)}\,q(\theta(t))
	= \sqrt{d_0}\big((1+s_0 t)\,q_0 + b_0 t\,w_0\big),
	\end{equation*}
	Therefore, denoting $v(t)=(1+s_0t)q_0+b_0tw_0$, we get
	\begin{align*}
		\Sigma(t)
		= d(t)q(\theta(t))q(\theta(t))^\top 
		= d_0 v(t) v(t)^\top 
	\end{align*}
	This agrees with the rank-one fixed-rank geodesics described in \cite{thanwerdas2023bures}, once one identifies the initial tangent vector via the differential of $\phi$ at $[\theta_0,d_0]$.
\end{example}

\begin{example}
	\label{ex:general-n-k1}
	We briefly indicate how Example~\ref{ex. n=2 k=1} generalizes to the case $n\ge2$, $k=1$. The associated-bundle structure is now
	\begin{equation*}
	M(n,1)=S^{n-1}\times_{\{\pm1\}}(0,+\infty),\qquad
	\pi_M:M(n,1)\to\R \P^{n-1},
	\end{equation*}
	with coordinates $([q,d],b,u)$ where $q\in S^{n-1}$, $d\in(0,+\infty)$, $b\in\R^{n-1}$, $u\in\R$. The metric becomes
	\begin{equation*}
	h_{[q,d]}(q_\perp b,u) = d\|b\|^2 + d s^2,
	\end{equation*}
	where $s$ is the scalar solving $2ds=u$.
	
	The geodesic system is
	\begin{equation*}
	\begin{cases}
		\dot{\mathbf{q}} = \mathbf{q}\,\mathbf{b}, & \mathbf{q}(0)=\mathbf{q}_0,\\[0.2em]
		\dot d = 2ds, & d(0)=d_0,\\[0.2em]
		\dot b = -2bs, & b(0)=b_0,\\[0.2em]
		\dot s = \|b\|^2 - s^2, & s(0)=s_0,
	\end{cases}
	\end{equation*}
	where
	\begin{equation*}
	\mathbf{q}=[q\;\;q_\perp]\in O(n),
	\qquad
	\mathbf{b}=\begin{bmatrix}
		0 & -b^\top \\
		b & 0
	\end{bmatrix}\in\mathfrak{o}(n).
	\end{equation*}
	As in Example~\ref{ex. n=2 k=1}, one shows that $\|b\|(r^2+1)$ is constant for $r=s/\|b\|$, which implies that $r$ evolves linearly in time and $b$, $s$, $d$ can be written explicitly. The evolution of $q$ is then obtained by integrating $\dot{\mathbf{q}}=\mathbf{q}\,\mathbf{b}$, which yields
	\begin{equation*}
	\mathbf{q}(t) = \mathbf{q}_0 \exp(\hat{\mathbf{b}}\theta(t)),
	\end{equation*}
	for a suitable angle function $\theta(t)$ and skew-symmetric matrix $\hat{\mathbf{b}}=\mathbf b/\|b\|$.
	The first column of this product gives $q(t)$ as a combination of $q_0$ and a fixed orthogonal direction, explicitly recovering the Bures--Wasserstein geodesics in the rank-one case $k=1$ for arbitrary $n$.
\end{example}

\subsubsection{Logarithm bijection}
\label{subsec. bijection}

Let $P_1,P_2\in\Gr(k,n)$, and choose $Q_1,Q_2\in\St(n,k)$ such that
\[
P_i=\Span(Q_i),\qquad i=1,2.
\]
The \emph{principal angles} between $P_1$ and $P_2$ are the uniquely determined numbers
\[
0\le \theta_1\le \cdots \le \theta_k\le \frac{\pi}{2}
\]
such that the singular values of $Q_1^\top Q_2$ are $\cos\theta_1,\dots,\cos\theta_k$
(see \cite{bjorck1973numerical}, in particular Theorem~1). We set
\begin{equation}\label{eq:r-def}
	r:=\#\{i:\theta_i=\pi/2\}.
\end{equation}
Thus, $r$ is the number of principal angles equal to $\pi/2$.

For the Grassmann manifold endowed with the Euclidean/Frobenius metric
(as in \cite[Section~5]{bendokat2024grassmann}), the set of minimizing logarithms
between $P_1$ and $P_2$ is indexed by
\begin{equation}\label{eq:I1}
	I_{\Gr}(P_1,P_2)
	:=
	\Bigl\{
	\diag(Q,I_{k-r}) : Q\in O(r)
	\Bigr\}.
\end{equation}
In particular, $I_{\Gr}(P_1,P_2)$ is naturally identified with $O(r)$.

Now let $\Sigma_1,\Sigma_2\in\Sym_+(n,k)$, and choose full column-rank factors
$X,Y\in\Mat(n,k)_*$ such that
\[
\Sigma_1=XX^\top,\qquad \Sigma_2=YY^\top.
\]
By \cite[Theorem~5.3]{thanwerdas2023bures}, the minimizing Bures--Wasserstein logarithms
between $\Sigma_1$ and $\Sigma_2$ are indexed by
\begin{equation}\label{eq:I2}
	I_{BW}(\Sigma_1,\Sigma_2)
	:=
	\bigl\{
	R\in O(k): X^\top YR\in\Cov(k)
	\bigr\}
	=
	\bigl\{
	R\in O(k): X^\top YR=(X^\top\Sigma_2X)^{1/2}
	\bigr\}.
\end{equation}
The next Theorem and Corollary show that the two index sets \eqref{eq:I1} and \eqref{eq:I2}
carry the same orthogonal degree of freedom.

\begin{theorem}
\label{thm:log-index-bijection}
Let $i\in\{1,2\}$. Let $\Sigma_i=Q_iD_iQ_i^\top \in \Sym^+(n,k)$ with
$Q_i\in\St(n,k)$ and $D_i\in\Sym_+(k)$, and set
$P_i:=\Imm(\Sigma_i)=\Span(Q_i)$. Let
\[
X:=Q_1D_1^{1/2},\qquad Y:=Q_2D_2^{1/2},
\]
so that $\Sigma_1=XX^\top$ and $\Sigma_2=YY^\top$, and define
\[
l:=\rk(X^\top Y).
\]
	Then:
	\begin{enumerate}
		\item[\textnormal{(i)}]one has
		\[
		k-l=r=\#\{i:\theta_i=\pi/2\},
		\]
		where $\theta_1,\dots,\theta_k$ are the principal angles between $P_1$ and $P_2$;
		\item[\textnormal{(ii)}] there is a bijection
		\[
		O(r)\ \longleftrightarrow\ I_{BW}(\Sigma_1,\Sigma_2).
		\]
	\end{enumerate}
\end{theorem}

\begin{proof}
	(i) Since $D_1^{1/2},D_2^{1/2}$ are invertible, we have
	\[
	\rk(X^\top Y)=\rk(D_1^{1/2}Q_1^\top Q_2D_2^{1/2})=\rk(Q_1^\top Q_2).
	\]
	The singular values of $Q_1^\top Q_2$ are $\cos\theta_i$, hence
	$\rk(Q_1^\top Q_2)=\#\{i:\cos\theta_i\neq 0\}=\#\{i:\theta_i\neq\pi/2\}$,
	so $k-l=r=\#\{i:\theta_i=\pi/2\}$.
	
	(ii) Let $X^\top Y=ULV^\top$ be a singular value decomposition with $U,V\in O(k)$ and
	\[
	L=\diag(\sigma_1,\dots,\sigma_l,0,\dots,0),\qquad \sigma_1,\dots,\sigma_l>0.
	\]
	Then 
	\[(X^\top \Sigma_2 X)^{1/2}=(X^\top YY^\top X)^{1/2}=ULU^\top.\]
	Let $R \in I_{BW}(\Sigma_1,\Sigma_2)$. The constraint $X^\top YR=(X^\top \Sigma_2 X)^{1/2}$ translates to
	\[
	ULV^\top R=ULU^\top,
	\]
	which, in view of the definition of $L$, is equivalent to
	\[V^\top R U=
	\begin{pmatrix}
		I_l & 0\\
		0 & R_r
	\end{pmatrix}\ \text{ with } R_r\in O(r)\]
	and therefore
	\begin{equation}\label{eq:R-param}
		R=V\begin{pmatrix}
			I_l & 0\\
			0 & R_r
		\end{pmatrix}U^\top,\qquad R_r\in O(r).
	\end{equation}
	Conversely, any $R$ of the form \eqref{eq:R-param} satisfies \eqref{eq:I2}. This yields the desired bijection.
\end{proof}

The next corollary follows from Theorem \ref{thm:log-index-bijection} and \eqref{eq:I1}.

\begin{corollary}[Bijection of index sets]
	\label{coroll:log-bij}
	Let $\Sigma_1, \Sigma_2 \in \Sym^+(n,k)$ and set $P_1=\Imm(\Sigma_1),\, P_2=\Imm(\Sigma_2)$.
	Then, there exists a bijection between the index sets $I_{\Gr}(P_1,P_2)$ and $I_{BW}(\Sigma_1,\Sigma_2)$.
\end{corollary}

\begin{remark}
	The bijection of Corollary \ref{coroll:log-bij} can be readily transferred from logarithms to length-minimizing geodesics on $\Sym^+(n,k)$ and $\Gr(k,n)$: each choice of $R_r\in O(r)$ corresponds simultaneously to a minimizing Grassmannian geodesic (via \cite{bendokat2024grassmann}) and to a minimizing Bures--Wasserstein geodesic on the fixed-rank covariance manifold (via \cite{thanwerdas2023bures}). This suggests a geometric connection between $(\Sym^+(n,k),d_{BW})$ and $(\Gr(k,n),d_E)$, where $d_E$ denotes the Euclidean distance.
	However, the horizontal part of the pullback metric on $M(n,k)$,
	\begin{equation*}
		h_{[Q,D]}^{\hor}(V,W) := \Tr\big(B_V D B_W^\top\big),
	\end{equation*}
	clearly differs from the Euclidean metric on $\Gr(k,n)$, as it carries vertical positional information through the anisotropic weight matrix $D$.
	
	This raises several questions and avenues for further investigation. How should one interpret this geodesic correspondence at the level of the metric structure? To what extent do Euclidean-geometric features of $\Gr(k,n)$ (such as distances, angles, curvature, cut times and loci) encode geometric or statistical properties of $(\Sym^+(n,k),d_{BW})$? More speculatively, one may ask whether the Bures--Wasserstein geometry of the fixed-rank stratum can be viewed as an anisotropically weighted perturbation of the Euclidean geometry of $\Gr(k,n)$, coupled with the Bures--Wasserstein metric itself on the full rank fibers; we leave these as open directions for future work.
\end{remark}

\section{Discussion and perspectives}

In this work we revisited the Bures--Wasserstein geometry of the fixed-rank stratum
\(\Sym^+(n,k)\) through the associated-bundle model
\begin{equation*}
\Sym^+(n,k) \;\cong\; M(n,k) \;=\; \St(n,k)\times_{O(k)} \Sym^{+}(k),
\end{equation*}
and we pulled back the Bures--Wasserstein metric to this bundle. This viewpoint makes
explicit the separation between motion of the image subspace, encoded by the Grassmannian
\(\Gr(k,n)\), and the covariance structure within that subspace, encoded by the SPD fiber
\(\Sym^{+}(k)\). Within this framework we derived an ODE system for Bures--Wasserstein geodesics in bundle coordinates, proved that the fibers \(M_P(n,k)\) are totally
geodesic,
and established a bijection between Bures--Wasserstein logarithms on \(\Sym^+(n,k)\) and
Grassmannian logarithms on \(\Gr(k,n)\).

From an information-geometric perspective, \(\Cov(n)\) can be identified with the space of
centred Gaussian measures on \(\mathbb{R}^n\), and the Bures--Wasserstein metric coincides
with the \(L^2\)-Wasserstein metric restricted to this Gaussian family. The fixed-rank strata
\(\Sym^+(n,k)\) then correspond to degenerate Gaussian distributions whose support lies in a
\(k\)-dimensional subspace. The bundle picture developed here emphasizes that the support
subspace and the covariance within that subspace are controlled by geometrically distinct
factors. We do not pursue concrete statistical models in this note, but we expect this
geometric separation to be a natural language for structured Gaussian families and
low-rank covariance estimators, such as in factor models or low rank covariance regularization. 

A first direction for further work concerns the horizontal--vertical splitting induced by the
Bures--Wasserstein metric. In the quotient description of \(\Sym^+(n,k)\) the vertical
directions correspond to deformations that preserve the image subspace, while the
horizontal directions encode motion of the subspace itself. In our associated-bundle model
this splitting becomes particularly transparent in terms of the variables \((Q,D,B,S)\). It
would be interesting to understand more systematically how this splitting interacts with
curvature and with the cut locus, and to what extent one can describe Bures--Wasserstein
geodesics by combining purely Grassmannian information with fiberwise covariance
evolution.

A related and more subtle issue is the link between Bures--Wasserstein geodesics on
\(\Sym^+(n,k)\) and minimizing geodesics on \(\Gr(k,n)\). We have shown that, for each choice
of index \(R_r \in O(r)\), there is a one-to-one correspondence between Grassmannian
logarithms (and minimizing geodesics) on \((\Gr(k,n),d_E)\) and Bures--Wasserstein
logarithms (and minimizing geodesics) on \((\Sym^+(n,k),d_{\mathrm{BW}})\). Here \(d_E\)
denotes the Euclidean metric on the Grassmannian induced by the ambient Frobenius
structure, which is substantially different from the horizontal part of the Bures--Wasserstein metric on \(\Sym^+(n,k)\). In particular, the projection
\(\pi_M : M(n,k)\to \Gr(k,n)\) is not isometric on horizontal spaces in the usual sense of
Riemannian submersions. A more conceptual explanation of why the index sets of
minimizing geodesics nevertheless coincide --- beyond the explicit matrix calculations used
here and in \cite{thanwerdas2023bures} --- is, to our knowledge, still missing. Clarifying
this relationship could shed further light on the role of \(\Gr(k,n)\) as a base space for the
Bures--Wasserstein geometry of low-rank covariance matrices.

Finally, the geodesic ODE system and the related conserved momentum in bundle coordinates suggest several numerical
questions. Expressing Bures--Wasserstein geodesics in terms of \((Q,D,B,S)\) makes it
possible in principle to construct integrators that combine existing algorithms on
\(\Gr(k,n)\) and \(\Sym^{+}(k)\) (for instance, Grassmannian exponential and logarithm
maps, or Sylvester solvers for \(S\)) into schemes on \(\Sym^+(n,k)\). Investigating the
numerical stability and efficiency of such schemes, and comparing them with methods
based directly on the quotient model \(\Mat(n,k)^*/O(k)\), is a natural continuation of the
present work. We leave these algorithmic developments for future research.

\appendix

\section{Appendix}
\label{sec:appendix}

\paragraph{Proof of Proposition \ref{prop. tangent spaces of associated bundle}.}

%
	
	Let $q:P\times F\to E$ be the canonical projection onto the associated bundle
	$E=P\times_G F$. By construction, $q$ is a submersion, so
	\[
	dq_{(p,f)}:T_pP\times T_fF\to T_{[p,f]}E
	\]
	is surjective and
	\[
	T_{[p,f]}E\cong \bigl(T_pP\times T_fF\bigr)\big/\Ker(dq_{(p,f)}).
	\]
	
	The kernel is the tangent space at $(p,f)$ to the $G$-orbit
	\[
	\mathcal O_{(p,f)}=\{(p,f)\cdot g:g\in G\}.
	\]
	For $\xi\in\mathfrak g$, consider the curve
	\begin{align*}
		\gamma_\xi(t)
		&:=(p,f)\cdot \exp(t\xi)\\
		&=\bigl(p\exp(t\xi),\,L_{\exp(-t\xi)}f\bigr),
	\end{align*}
	which lies in $\mathcal O_{(p,f)}$. Differentiating at $t=0$ gives
	\[
	\dot\gamma_\xi(0)=(p\xi,-L_\xi^*f).
	\]
	Since every tangent vector to $\mathcal O_{(p,f)}$ at $(p,f)$ arises in this way, we obtain
	\[
	\Ker(dq_{(p,f)})
	=
	T_{(p,f)}\mathcal O_{(p,f)}
	=
	\{(p\xi,-L_\xi^*f):\xi\in\mathfrak g\}.
	\]
	Therefore
	\[
	T_{[p,f]}E
	=
	\bigl(T_pP\times T_fF\bigr)\Big/\{(p\xi,-L_\xi^*f):\xi\in\mathfrak g\},
	\]
	as claimed.

\paragraph{Proof of Lemma \ref{lemma:local to global diff}.}

	We first show that $\phi$ is a bijection.  
	
	\emph{Surjectivity.} Let $y\in N$ and set $b := \pi_N(y)\in B$. The restriction
	\[
	\phi_b := \phi\big|_{\pi_M^{-1}(b)} : \pi_M^{-1}(b) \to \pi_N^{-1}(b)
	\]
	is a diffeomorphism by assumption. In particular, $y\in \pi_N^{-1}(b)$ admits a preimage $x\in \pi_M^{-1}(b)$ with $\phi(x)=y$.
	
	\emph{Injectivity.} Suppose $\phi(x_1)=\phi(x_2)$ for $x_1,x_2\in M$. Then
	\begin{equation*}
	\pi_M(x_1) = \pi_N(\phi(x_1)) = \pi_N(\phi(x_2)) = \pi_M(x_2),
	\end{equation*}
	so $x_1$ and $x_2$ lie in the same fiber of $\pi_M$. Since the restriction of $\phi$ to that fiber is a diffeomorphism, it is injective there, and thus $x_1=x_2$.
	
	We have shown that $\phi$ is a bijection.
	
	Since $M$ and $N$ are smooth fiber bundles and their fibers are diffeomorphic to each other, there exists an open cover $\{U_\alpha\}$ of $B$ and smooth local trivializations
	\begin{equation*}
	\psi_\alpha^M: \pi_M^{-1} (U_\alpha) \to U_\alpha \times F,
	\qquad
	\psi_\alpha^N: \pi_N^{-1} (U_\alpha) \to U_\alpha \times F,
	\end{equation*}
	where $F$ is a model fiber diffeomorphic to the fibers of both bundles. In these coordinates we can write
	\begin{equation*}
	\widetilde{\phi}_\alpha
	:= \psi_\alpha^N \circ \phi\circ (\psi_\alpha^M)^{-1}: U_\alpha \times F \to U_\alpha \times F.
	\end{equation*}
	Because $\phi$ is a bundle morphism, i.e., it respects fibers, $\widetilde{\phi}_\alpha$ takes the form
	\begin{equation*}
	\widetilde{\phi}_\alpha(b,x)=(b,f_\alpha(b,x)),
	\end{equation*}
	where $f_\alpha: U_\alpha \times F \to F$ is smooth and, for each $b\in U_\alpha$, the map
	\begin{equation*}
	f_{\alpha,b} := f_\alpha(b,\cdot) : F\to F
	\end{equation*}
	is a diffeomorphism (since $\phi$ restricts to a diffeomorphism on each fiber).
	
	The differential of $\widetilde{\phi}_\alpha$ at $(b,x)$ has block form
	\begin{equation*}
	D\widetilde{\phi}_\alpha(b,x)
	= \begin{bmatrix}
		I & 0\\
		\ast & D_x f_{\alpha,b}
	\end{bmatrix},
	\end{equation*}
	where $I$ is the identity on $T_b U_\alpha$ and $D_x f_{\alpha,b}$ is the differential of $f_{\alpha,b}$ at $x$. Because each $f_{\alpha,b}$ is a diffeomorphism, $D_x f_{\alpha,b}$ is invertible, so $D\widetilde{\phi}_\alpha(b,x)$ is invertible as well. Thus $\widetilde{\phi}_\alpha$ is a local diffeomorphism, and since it is also globally bijective on $U_\alpha\times F$, it is a global diffeomorphism there, with inverse
	\begin{equation*}
	\widetilde{\phi}_\alpha^{-1}(b,x) = (b,f_{\alpha,b}^{-1}(x)).
	\end{equation*}
	
	It follows that $\phi$ is smooth and has a smooth inverse in each local trivialization chart, and hence globally on $M$ and $N$. Therefore $\phi$ is a diffeomorphism.

\paragraph{Proof of Proposition \ref{prop:principal bundlle submersion}.}

	Let $Q\in\St(n,k)$, and let $Q_\perp$ be a completion to an orthonormal basis of $\R^n$, so that $[Q\; Q_\perp]\in O(n)$. Consider a horizontal tangent vector $V\in T_Q\St(n,k)$, which by definition can be written as $V=Q_\perp B$ for some $B\in\Mat(n-k,k)$.
	
	Let $Q(t)\in\St(n,k)$ be a smooth curve with $Q(0)=Q$ and $\dot Q(0)=V$. For any $x\in\R^k$, consider the curve $v(t)=Q(t)x\in\R^n$, which lies in the evolving subspace $P(t):=\pi_G(Q(t))\in\Gr(k,n)$. Then
	\begin{equation*}
	\dot v(0) = \dot Q(0)x = Q_\perp B x.
	\end{equation*}
	Since $x\in\R^k$ is arbitrary, the derivative $d\pi_G(Q)(V)$, viewed in the coordinates induced by the orthonormal basis $[Q\; Q_\perp]$ of $\R^n$, is exactly the matrix $B\in\Mat(n-k,k)$. In particular,
	\begin{equation*}
	d\pi_G(Q) : \cl{H}_Q \to T_{\Span(Q)}\Gr(k,n)\cong \Mat(n-k,k)
	\end{equation*}
	is an isomorphism, and for $V_1=Q_\perp B_1$ and $V_2=Q_\perp B_2$ we have
	\begin{equation*}
	g^{\St}_Q(V_1,V_2) = \Tr(B_1 B_2^\top ) = g^{\Gr}_{\Span(Q)}\big(d\pi_G(Q)(V_1),d\pi_G(Q)(V_2)\big),
	\end{equation*}
	by the definitions in \eqref{metric on grassmannian}. Thus $d\pi_G(Q)$ restricts to an isometry from $\cl{H}_Q$ onto $T_{\Span(Q)}\Gr(k,n)$, so $\pi_G$ is a Riemannian submersion.

\paragraph{Proof of Proposition \ref{prop:pullback-metric}.}

	Take a smooth curve $[Q(t),D(t)] \in M(n,k)$ such that
	\begin{equation*}
	Q(0)=Q,\quad \dot Q(0)=Q_\perp B,
	\qquad
	D(0)=D,\quad \dot D(0)=T,
	\end{equation*}
	with $B \in \Mat(n-k,k)$ and $T \in \Sym(k)$. Then the differential $d\phi_{[Q,D]} \colon T_{[Q,D]}M(n,k) \longrightarrow T_{QDQ^\top }\Sym^+(n,k)$ is given by the equation
	\begin{align*}
		d\phi_{[Q,D]}(Q_\perp B,T)
		&= \frac{d}{dt}\Big|_{t=0} Q(t)D(t)Q(t)^\top  \\
		&= Q_\perp B D Q^\top  + Q D B^\top  Q_\perp^\top  + Q T Q^\top .
	\end{align*}
	Denote
	\begin{equation*}
	\widetilde{V} := d\phi_{[Q,D]}(V),\qquad \widetilde{W} := d\phi_{[Q,D]}(W).
	\end{equation*}
	We are going to analyse separately the horizontal $g_\Sigma^{\hor}(\widetilde V,\widetilde W)=\Tr\big(P_\Sigma^\perp \widetilde V \Sigma^\dagger \widetilde W\big)$ and vertical $g_\Sigma^{\ver}(\widetilde V,\widetilde W)
	= \Tr\left(S_{\Sigma,\widetilde V}\,\Sigma\,S_{\Sigma,\widetilde W}\right)$ parts of the metric.

	\emph{Horizontal part.}
	Using $P_\Sigma^\perp=Q_\perp Q_\perp^\top $ and $\Sigma^\dagger = Q D^{-1} Q^\top $, we compute
	\begin{align*}
		P_\Sigma^\perp \widetilde{V}
		&= Q_\perp Q_\perp^\top 
		\big(Q_\perp B_V D Q^\top  + Q D B_V^\top  Q_\perp^\top  + Q T_V Q^\top \big) \\
		&= Q_\perp B_V D Q^\top ,
	\end{align*}
	since $Q_\perp^\top  Q = 0$ and $Q_\perp^\top  Q_\perp = I_{n-k}$.
	Thus
	\begin{equation*}
	P_\Sigma^\perp \widetilde{V}\,\Sigma^\dagger
	= Q_\perp B_V D Q^\top  Q D^{-1} Q^\top 
	= Q_\perp B_V Q^\top .
	\end{equation*}
	Multiplying by $\widetilde{W}$ gives
	\begin{align*}
		P_\Sigma^\perp \widetilde{V}\,\Sigma^\dagger \widetilde{W}
		&= Q_\perp B_V Q^\top 
		\big(Q_\perp B_W D Q^\top  + Q D B_W^\top  Q_\perp^\top  + Q T_W Q^\top \big) \\
		&= Q_\perp B_V D B_W^\top  Q_\perp^\top  + Q_\perp B_V T_W Q^\top .
	\end{align*}
	Taking the trace,
	\begin{align*}
		g_\Sigma^{\hor}(\widetilde{V},\widetilde{W})
		&= \Tr\big(P_\Sigma^\perp \widetilde{V}\,\Sigma^\dagger \widetilde{W}\big) \\
		&= \Tr\big(Q_\perp B_V D B_W^\top  Q_\perp^\top \big)
		+ \Tr\big(Q_\perp B_V T_W Q^\top \big).
	\end{align*}
	The second term vanishes thanks to cyclicity of the trace:
	\begin{equation*}
	\Tr\big(Q_\perp B_V T_W Q^\top \big)
	= \Tr\big(Q^\top  Q_\perp B_V T_W\big)
	= \Tr\big(0\cdot B_V T_W\big)
	= 0,
	\end{equation*}
	so
	\begin{equation*}
	g_\Sigma^{\hor}(\widetilde{V},\widetilde{W})
	= \Tr\big(Q_\perp B_V D B_W^\top  Q_\perp^\top \big)
	= \Tr\big(B_V D B_W^\top \big),
	\end{equation*}
	again by cyclicity of the trace.

	\emph{Vertical part.}
	For the vertical component, we use
	\begin{equation*}
	\widetilde{V}
	= Q_\perp B_V D Q^\top  + Q D B_V^\top  Q_\perp^\top  + Q T_V Q^\top ,
	\end{equation*}
	and obtain
	\begin{equation*}
	Q^\top  \widetilde{V} Q
	= Q^\top  Q_\perp B_V D Q^\top  Q 
	+ Q^\top  Q D B_V^\top  Q_\perp^\top  Q
	+ Q^\top  Q T_V Q^\top  Q
	= T_V,
	\end{equation*}
	since $Q^\top  Q = I_k$ and $Q^\top  Q_\perp=0$.
	Thus
	\begin{equation*}
	S_{\Sigma,\widetilde{V}} = Q \cl S_D(T_V) Q^\top ,
	\end{equation*}
	and similarly $S_{\Sigma,\widetilde{W}} = Q \cl S_D(T_W) Q^\top $. Therefore,
	\begin{align*}
		g_\Sigma^{\ver}(\widetilde{V},\widetilde{W})
		&= \Tr\big(S_{\Sigma,\widetilde{V}}\,\Sigma\,S_{\Sigma,\widetilde{W}}\big) \\
		&= \Tr\big(Q \cl S_D(T_V) Q^\top \, Q D Q^\top \, Q \cl S_D(T_W) Q^\top \big) \\
		&= \Tr\big(\cl S_D(T_V) D \cl S_D(T_W)\big),
	\end{align*}
	again by cyclicity of the trace.
	
	Putting the horizontal and vertical contributions together, we obtain
	\begin{equation*}
	h_{[Q,D]}(V,W)
	= g_\Sigma^{BW(n,k)}(\widetilde{V},\widetilde{W})
	= \Tr\big(B_V D B_W^\top \big)
	+ \Tr\big(\cl S_D(T_V)\,D\,\cl S_D(T_W)\big),
	\end{equation*}
	which proves \eqref{eq:pullback-metric}.

\paragraph{Proof of Theorem \ref{thm:BW-geodesic-ODE}.}

	With the metric $h$ from Subsection~\ref{subsec. metric on sym}, the kinetic energy (Lagrangian) is
	\begin{align*}
		L(B,D,T)
		&:= \frac{1}{2}h\big(\dot\gamma(t),\dot\gamma(t)\big) \\
		&= \frac{1}{2}\Tr\big(B D B^\top \big) + \frac{1}{2}\Tr\big(\cl S_D(T)\,D\,\cl S_D(T)\big) \\
		&=: L^{\hor}(B,D) + L^{\ver}(D,T),
	\end{align*}
	where $\cl S_D(T)=: S$ is the unique solution of the Sylvester equation
	\begin{equation*}
	D S + S D = T.
	\end{equation*}
	
	A geodesic is a smooth curve $\gamma$ such that $L$ is stationary under all smooth variations with fixed endpoints. Writing $L$ in terms of the coordinates $(B,D,T)$ and using that $L$ does not depend on $Q$, the Euler--Lagrange equations reduce to
	\begin{equation}
		\label{eq:EL-basic}
		\frac{d}{dt}\big(\partial_B L\big)=0,
		\qquad
		\frac{d}{dt}\big(\partial_T L\big) - \partial_D L=0.
	\end{equation}
	
	We now compute the partial derivatives. Throughout we use the straightforward matrix-variational identity
	\begin{equation*}
	\delta f(N)=\Tr\big(\delta N\,(\partial_N f)^\top \big),
	\end{equation*}
	for a scalar function $f$ of a matrix argument $N$, where $\delta$ denotes first variation; see \cite{magnus2019matrix} for a general reference on matrix calculus.
	
	\emph{Derivative with respect to $B$.}
	We have
	\begin{align*}
		\delta_B \Tr(B D B^\top )
		&= \frac{d}{d\varepsilon}\Tr\big((B+\varepsilon\,\delta B)D(B+\varepsilon\,\delta B)^\top \big)\big|_{\varepsilon=0} \\
		&= \Tr(\delta B D B^\top ) + \Tr(B D \delta B^\top ) \\
		&= 2\Tr(\delta B D B^\top ),
	\end{align*}
	so
	\begin{equation*}
	\partial_B L = \partial_B L^{\hor} = B D.
	\end{equation*}
	The first Euler--Lagrange equation in \eqref{eq:EL-basic} therefore gives
	\begin{equation*}
		\label{eq:BD-constant}
		\frac{d}{dt}(BD)=0
		\qquad\Longrightarrow\qquad
		BD = K,
	\end{equation*}
	for some constant matrix $K\in\Mat(n-k,k)$ along the geodesic. Since $D$ is invertible on $\Sym^+(k)$, we may write
	\begin{equation}
		\label{eq:Bdot}
		\dot B = - B\,\dot D\,D^{-1} = - B T D^{-1},
	\end{equation}
	using the relation $\dot D=T$ (the vertical component $T$ is precisely the derivative of $D$ along $\gamma$). In terms of the Sylvester variable $S=\cl S_D(T)$ we can equivalently write
	\begin{equation*}
	\dot B = -B(DS+SD)D^{-1}.
	\end{equation*}
	
	\emph{Derivative with respect to $T$.}
	For the vertical part we compute
	\begin{align*}
		\delta_T \Tr\big(\cl S_D(T)\,D\,\cl S_D(T)\big)
		&= \Tr\big(\cl S_D(\delta T)\,D\,\cl S_D(T)\big)
		+ \Tr\big(\cl S_D(T)\,D\,\cl S_D(\delta T)\big) \\
		&= \Tr\big(\cl S_D(\delta T)(D \cl S_D(T)+\cl S_D(T)D)\big) \\
		&= \Tr\big(\delta T\, \cl S_D(D \cl S_D(T)+\cl S_D(T)D)\big)\\
		&= \Tr\big(\delta T\, \cl S_D(T)\big),
	\end{align*}
	where we used the linearity and self-adjointness of $\cl S_D(\cdot)$ with respect to the Frobenius inner product, and the Sylvester equation
	\begin{equation*}
	D \cl S_D(T)+\cl S_D(T)D = T.
	\end{equation*}
	For an in-depth treatment of the Sylvester and Lyapunov equations, where all the corresponding properties are derived, see \cite{bhatia2007positive}.
	Thus
	\begin{equation*}
	\partial_T L = \partial_T L^{\ver} = \frac{1}{2}\cl S_D(T).
	\end{equation*}
	
	\emph{Derivative with respect to $D$.}
	For the horizontal part,
	\begin{equation*}
	\delta_D \Tr(B D B^\top ) = \Tr(B^\top  B\,\delta D),
	\end{equation*}
	so
	\begin{equation*}
	\partial_D L^{\hor} = \frac{1}{2}B^\top  B.
	\end{equation*}
	For the vertical part we differentiate the Sylvester constraint
	\begin{equation*}
	D \cl S_D(T) + \cl S_D(T) D = T
	\end{equation*}
	at fixed $T$, obtaining
	\begin{equation*}
	D\,\delta_D \cl S_D(T) + \delta_D \cl S_D(T) D
	= -\delta D\,\cl S_D(T) - \cl S_D(T)\,\delta D,
	\end{equation*}
	and hence
	\begin{equation*}
	\delta_D \cl S_D(T)
	= -\cl S_D\big(\delta D\,\cl S_D(T) + \cl S_D(T)\,\delta D\big).
	\end{equation*}
	A direct calculation then gives
	\begin{align*}
		\delta_D \Tr\big(\cl S_D(T)\,D\,\cl S_D(T)\big)
		&= \Tr\big(\delta_D \cl S_D(T)\,D\,\cl S_D(T)\big)
		+ \Tr\big(\cl S_D(T)\,\delta D\,\cl S_D(T)\big) \\
		&\quad + \Tr\big(\cl S_D(T)\,D\,\delta_D \cl S_D(T)\big) \\
		&= \Tr\big(\cl S_D(T)\,\delta D\,\cl S_D(T)\big)
		- \Tr\big(\delta D\,\cl S_D(T)^2\big) \\
		&= -\Tr\big(\delta D\,\cl S_D(T)^2\big),
	\end{align*}
	so
	\begin{equation*}
	\partial_D L^{\ver} = -\frac{1}{2}\cl S_D(T)^2.
	\end{equation*}
	Combining the two contributions,
	\begin{equation*}
	\partial_D L = \frac{1}{2}\big(B^\top  B - \cl S_D(T)^2\big).
	\end{equation*}
	
	\emph{Evolution equation for $S=\cl S_D(T)$.}
	Plugging the above expressions into the second Euler--Lagrange equation in \eqref{eq:EL-basic},
	\begin{equation*}
	\frac{d}{dt}(\partial_T L) - \partial_D L = 0,
	\end{equation*}
	we obtain
	\begin{equation*}
	\frac{1}{2}\frac{d}{dt}\cl S_D(T) - \frac{1}{2}\big(B^\top  B - \cl S_D(T)^2\big) = 0,
	\end{equation*}
	or, writing $S=\cl S_D(T)$,
	\begin{equation}
		\label{eq:Sdot}
		\dot S = B^\top  B - S^2.
	\end{equation}
	Finally, the evolution of $\mathbf{Q}$ is given by the horizontal component $Q_\perp B$ of the velocity and the skew-symmetric matrix $\mathbf{B}$ as in the statement. Conversely, any solution of \eqref{eq:full-geodesic-system} arises from a critical point of the energy functional with fixed endpoints, hence from a geodesic. This establishes the equivalence.
	As a final note, we observe that $O(k)$-invariance of the system naturally follows from the $O(k)$-invariance of the metric.

\printbibliography

\end{document}